\documentclass[11pt]{article}
\usepackage{latexsym}

\newtheorem{theorem}{Theorem}[section]
\newtheorem{lemma}[theorem]{Lemma}
\newtheorem{proposition}[theorem]{Proposition}
\newtheorem{definition}{Definition}
\newtheorem{corollary}[theorem]{Corollary}

\newcommand{\qed}{\ \hfill\mbox{$\Box$}\vspace{\baselineskip}}
\newenvironment{proof}{\noindent {\bf Proof:}}{{\qed}}

\newcommand{\conv}{\mbox{{\rm conv}}}
\begin{document}

\title{Shelling and triangulating\\ the (extra)ordinary polytope}

\author{Margaret M. Bayer\thanks{This research was supported by the sabbatical
        leave program of the University of Kansas, and was conducted while the 
        author was at the Mathematical Sciences Research Institute, supported in
        part by NSF grant DMS-9810361, and at 
        Technische Universit\"{a}t Berlin, supported in part by 
        Deutsche Forschungs-Gemeinschaft, through the
        DFG Research Center ``Mathematics for Key Technologies'' (FZT86) and
       the Research Group ``Algorithms, Structure, Randomness'' (FOR~13/1-1).}\\
        Department of Mathematics\\
        University of Kansas\\
        Lawrence, KS  66045-7523\\
        bayer@math.ukans.edu}

\date{}

\maketitle

\begin{abstract}
Ordinary polytopes were introduced by Bisztriczky as a (nonsimplicial)
generalization of cyclic polytopes.
We show that the colex order of facets of the ordinary polytope is a
shelling order.  This shelling shares many nice properties with the
shellings of simplicial polytopes.  
We also give a shallow triangulation of the ordinary polytope, and show
how the shelling and the triangulation are used to compute the
toric $h$-vector of the ordinary polytope.
As one consequence, we get that the contribution from each shelling component
to the $h$-vector is nonnegative.
Another consequence is a combinatorial proof that the entries of the 
$h$-vector of any ordinary polytope are simple sums of binomial coefficients.
\end{abstract}

\section{Introduction}
\subsection{Motivation}
This paper has a couple of main motivations.
The first comes from the study of toric $h$-vectors of convex polytopes.
The $h$-vector played a crucial role in the characterization of face vectors
of simplicial polytopes \cite{billera-lee,McMullen,sta80}.
In the simplicial case, the $h$-vector is
linearly equivalent to the face vector, and has a
combinatorial interpretation in a shelling of the polytope.
The $h$-vector of a simplicial polytope is also the sequence of
Betti numbers of the associated toric variety.
In this context it generalizes to nonsimplicial polytopes. 
However, for nonsimplicial polytopes, we do not have a good combinatorial
understanding of the entries of the $h$-vector.

The definition of the (toric) $h$-vector for general polytopes
(and even  more generally, for Eulerian posets) first appeared in
\cite{sta87}.
Already there Stanley raised the issue of computing the $h$-vector from a
shelling of the polytope.
Associated with any shelling, $F_1$, $F_2$, \ldots, $F_n$,
of a polytope $P$ is a partition of the
faces of $P$ into the sets ${\cal G}_j$ of faces of $F_j$ not in
$\cup_{i<j} F_i$.
The $h$-vector can be decomposed into contributions from each set
${\cal G}_j$.
When $P$ is simplicial, the set ${\cal G}_j$ is a single interval
$[G_j,F_j]$ in the face lattice of $P$, and the contribution to
the $h$-vector is a single 1 in position $|G_j|$.
For nonsimplicial polytopes, the set ${\cal G}_j$ is not so simple.
It is not clear whether the contribution to the $h$-vector 
from ${\cal G}_j$ must be nonnegative, and, if it is, whether
it counts something natural.
(Tom Braden \cite{braden} has announced a positive answer to this question,
based on \cite{Barthel,Karu}.)
Another issue is the relation of the $h$-vector of a polytope $P$ to the
$h$-vector of a triangulation of $P$.
This is addressed in \cite{Bayer-weakly,sta92}.

A problem in studying nonsimplicial polytopes is the difficulty of generating
examples with a broad range of combinatorial types.
Bisztriczky \cite{Bisz} discovered the fascinating ``ordinary''
polytopes, a class of generally nonsimplicial polytopes, which includes
as its simplicial members the cyclic polytopes.
These polytopes have been studied further in 
\cite{bayer-ordflags,Bayer-Bru-Ste,Dinh}.
In particular, in \cite{bayer-ordflags}, it is shown that ordinary
polytopes have surprisingly nice $h$-vectors, namely, the $h$-vector
is the sum of the $h$-vector of a cyclic polytope and the shifted
$h$-vector of a lower-dimensional cyclic polytope.
These $h$-vectors were calculated from the flag vectors, and the
calculation did not give a combinatorial explanation for the nice
form that came out.
So we were motivated to find a combinatorial interpretation for
these $h$-vectors, most likely through shellings or triangulations of
the polytopes.

This paper is organized as follows.
In the second part of this introduction we give the main definitions.
The brief Section~2 warms up with the natural triangulation of the multiplex.
Section~3 is devoted to showing that the colex order of facets is a shelling
of the ordinary polytope.
The proof, while laborious, is constructive, explicitly describing the
minimal new faces of the polytope as each facet is shelled on.
We then turn in Section~4 to $h$-vectors of multiplicial polytopes in
general, and of the ordinary polytope in particular.
Here a ``fake simplicial $h$-vector'' arises in the shelling of the
ordinary polytope.
In Section~5, the triangulation of the multiplex is used to triangulate
the boundary of the ordinary polytope.
This triangulation is shown to have a shelling compatible with the
shelling of Section~3.
The shelling and triangulation together
explain combinatorially the $h$-vector of the ordinary polytope.

Finally, a comment about the title of this paper.
Bisztriczky named these polytopes ``ordinary polytopes'' to invoke
the idea of ordinary curves.
The name is, of course, a bit misleading, as it is applied to a
truly extraordinary class of polytopes.
We feel that these polytopes are extraordinary
because of their special structure, but we hope that they will also
turn out to be extraordinary for their usefulness in understanding
general convex polytopes.

\subsection{Definitions}
For common polytope terminology, refer to \cite{Ziegler}.

The {\em toric $h$-vector} was defined by Stanley for
Eulerian posets, including the face lattices of convex polytopes.
\begin{definition}[\cite{sta87}] {\em
Let $P$ be a $(d-1)$-dimensional polytopal sphere.
The $h$-vector and $g$-vector of $P$ are encoded as polynomials:
$h(P,x)=\sum_{i=0}^d h_ix^{d-i}$ and
$g(P,x)=\sum_{i=0}^{\lfloor d/2\rfloor}g_ix^i$, with the relations $g_0=h_0$
and $g_i=h_i-h_{i-1}$ for $1\le i\le d/2$. 
Then the $h$-polynomial and $g$-polynomial are defined by the recursion
\begin{enumerate}
\item $g(\emptyset,x)=h(\emptyset,x)=1$, and
\item $\displaystyle h(P,x)=\sum_{{\mbox{{\scriptsize $G$ face of $P$}}}\atop
       {G\not=P}} g(G,x)(x-1)^{d-1-\dim G}$.
\end{enumerate}
}\end{definition}

It is easy to see that the $h$-vector depends linearly on the flag vector.
In the case of simplicial polytopes, the formulas reduce to the well-known
transformation between $f$-vector and $h$-vector.

\begin{definition}[\cite{Ziegler}]{\em
Let $\cal C$ be a pure $d$-dimensional polytopal complex.
If $d=0$, then a {\em shelling} of $\cal C$ is any ordering
of the points of $C$.
If $d>0$, then
a {\em shelling} of $\cal C$ is a linear ordering $F_1$, $F_2$, \ldots,
$F_s$ of the facets of $\cal C$ such that for $2\le j\le s$, $F_j\cap
(\cup_{i<j}F_i)$ is nonempty and is the union of ridges of $\cal C$ that
form the initial segment of a shelling of $F_j$.
}\end{definition}

\begin{definition}[\cite{Bayer-weakly}]
{\em A triangulation $\Delta$ of a polytopal complex $\cal C$
is {\em shallow} if and only if every face $\sigma$ of $\Delta$ is contained in
a face of $\cal C$ of dimension at most $2\dim\sigma$.
}\end{definition}

\begin{theorem}[\cite{Bayer-weakly}]
If $\Delta$ is a simplicial sphere forming a shallow triangulation of the 
boundary of the convex $d$-polytope $P$, then $h(\Delta,x)=h(P,x)$.
\end{theorem}
Note:
in \cite{Bayer-weakly} Theorem~4 gives $h(P,x)=h(\Delta,x)$ for a shallow
subdivision $\Delta$ of the solid polytope $P$.
The proof goes through for shallow subdivisions of the boundary, because it
is based on the uniqueness of low-degree acceptable functions \cite{sta87}, 
which holds for lower Eulerian posets.

\begin{definition}[\cite{Bisz-mult}] \label{def-multi} {\em
A $d$-dimensional {\em multiplex} is a polytope with an ordered list of
vertices,
$x_0$, $x_1$, \ldots, $x_n$, with facets $F_0$, $F_1$, \ldots, $F_n$ given by
$$F_i = \conv\{x_{i-d+1}, x_{i-d+2}, \ldots, x_{i-1}, x_{i+1}, x_{i+2}, \ldots,
x_{i+d-1}\},$$
with the conventions that $x_i=x_0$ if $i<0$, and $x_i=x_n$ if $i>n$.
}\end{definition}

Given an ordered set $V=\{x_0,x_1,\ldots,x_n\}$, a subset $Y\subseteq V$ is
called a {\em Gale subset} if between any two elements of $V\setminus Y$
there is an even number of elements of $Y$.
A polytope $P$ with ordered vertex set $V$ is a {\em Gale polytope}
if the set of vertices of each facet is a Gale subset.
\begin{definition}[\cite{Bisz}] {\em
An {\em ordinary polytope} is a Gale polytope such that each facet is a
multiplex with the induced order on the vertices.
}\end{definition}
Cyclic polytopes can be characterized as the simplicial Gale polytopes.
Thus the only simplicial ordinary polytopes are cyclics.
In fact, these are the only ordinary polytopes in even dimensions.
However, the odd-dimensional, nonsimplicial ordinary polytopes are quite 
interesting.

We use the following notational conventions.
Vertices are generally denoted by integers $i$ rather than by $x_i$.
Where it does not cause confusion, a face of a polytope or a 
triangulation is identified with its vertex set.
Interval notation is used to denote sets of consecutive
integers, $[a,b]=\{a,a+1,\ldots, b-1,b\}$.

\section{Triangulating the multiplex}
\label{triang-multi}
Multiplexes have minimal triangulations that are particularly easy to
describe.
\begin{theorem}\label{trian-multi}
Let $M^{d,n}$ be a multiplex with ordered vertices 
0, 1, \ldots, $n$.
For $0\le i\le n-d$, let $T_i$ be the convex hull of $[i,i+d]$.
Then $M^{d,n}$ has a shallow triangulation as the union of the $n-d+1$
$d$-simplices $T_i$.
\end{theorem}
\begin{proof}
The proof is by induction on $n$.
For $n=d$, the multiplex $M^{d,d}$ is the simplex $T_0$ itself.
Assume $M^{d,n}$ has triangulation into simplices $T_i$, $0\le i\le n-d$.
Consider the multiplex $M^{d,n+1}$ with ordered vertices 
0, 1, \ldots, $n$.
Then $M^{d,n+1}=\conv(M^{d,n}\cup\{n+1\})$, where $n+1$ is a point
beyond facet $F_n$ of $M^{d,n}$, beneath the facets $F_i$ for $0\le i\le n-d+1$,
and in the affine hulls of the facets $F_i$ for $n-d+2\le i\le n-1$.
(See \cite{Bisz-mult}.)
Thus, $M^{d,n+1}$ is the union of $M^{d,n}$ and 
$\conv(F_n\cup\{n+1\})=T_{n+1-d}$, and $M^{d,n}\cap T_{n+1-d}=F_n$.
By the induction assumption, 
the simplices $T_i$, with $0\le i\le n+1-d$, form a 
triangulation of $M^{d,n+1}$.

The dual graph of the triangulation is simply a path. 
(The dual graph is the graph having a vertex for each
$d$-simplex, and an edge between two vertices if the corresponding
$d$-simplices share a $(d-1)$-face.) 
The ordering $T_0$, $T_1$, $T_2$, \ldots, $T_{n-d}$ is a shelling of the
simplicial complex that triangulates $M^{d,n}$.
So the $h$-vector of the triangulation is $(1, n-d, 0, 0, \ldots)$.
This is the same as the $g$-vector of the boundary of the multiplex,  
which is the $h$-vector of the solid multiplex.
So by \cite{Bayer-weakly}, the triangulation is shallow.
\end{proof}

Note, however, that for $n\ge d+2$, $M^{d,n}$ is not {\em weakly neighborly}
(as observed in \cite{Bayer-Bru-Ste}). 
This means that it has nonshallow triangulations.
This is easy to see because the vertices 0 and $n$ are
not contained in a common proper face of $M^{d,n}$.

Consider the induced triangulation of the boundary of $M^{d,n}$.
For notational purposes we consider $T_0$ and $T_n$ separately.
All facets of $T_0$ except $[1,d]$ are boundary facets of $M^{d,n}$.
Write $T_{0\setminus 0}=[0,d-1]=F_0$, and $T_{0\setminus j}=[0,d]\setminus\{j\}$
for $1\le j\le d-1$.
Write $T_{n-d\setminus n}=[n-d+1,n]=F_n$, and
$T_{n\setminus j}=[n-d,n]\setminus\{j\}$
for $n-d+1\le j\le n-1$.
For $1\le i\le n-d-1$, the facets of $T_i$ are 
$T_{i\setminus j}=[i,i+d]\setminus\{j\}$.
Two of these facets ($j=i$ and $j=i+d$) intersect the interior of $M^{d,n}$.
For $1\le j\le n-1$,
the facet $F_j$ is triangulated by $T_{i\setminus j}$ for $j-d+1\le i\le j-1$
(and $0\le i \le n-d$).
The facet order $F_0$, $F_1$, \ldots, $F_n$, is a shelling of 
the multiplex $M^{d,n}$.
The $(d-1)$-simplices $T_{i\setminus j}$ in the 
order $T_{0\setminus 0}$, $T_{0\setminus 1}$, $T_{0\setminus 2}$,
$T_{1\setminus 2}$, \ldots, $T_{n-d-1\setminus n-2}$, $T_{n-d\setminus n-2}$, 
$T_{n-d\setminus n-1}$, $T_{n-d+1\setminus n}$
(increasing order of $j$ and, for each $j$, increasing order of $i$),
form a shelling of the triangulated boundary of $M^{d,n}$.

\section{Shelling the ordinary polytope}

Shelling is used to calculate the $h$-vector, and hence the $f$-vector
of simplicial complexes (in particular, the boundaries of simplicial polytopes).
This is possible because (1) the $h$-vector has a simple expression in terms
of the $f$-vector and vice versa; (2) in a shelling of a simplicial complex,
among the faces added to the subcomplex as a new facet is shelled on,
there is a unique minimal face; (3) the interval from this minimal new face 
to the facet is a Boolean algebra; and (4) the numbers of new faces given by
(3) match the coefficients in the $f$-vector/$h$-vector formula.
These conditions all fail for shellings of arbitrary polytopes.
However, some hold for certain shellings of ordinary polytopes.

As mentioned earlier, noncyclic ordinary polytopes exist only in odd
dimensions.
Furthermore, three-dimensional ordinary polytopes are quite different combinatorially from those in higher dimensions.
We thus restrict our attention to ordinary polytopes of odd dimension at least 
five.
It turns out that these are classified by the vertex figure of the first vertex.
\begin{theorem}[\cite{Bisz,Dinh}]
For each choice of integers $n\ge k\ge d=2m+1\ge 5$, there is a unique 
combinatorial type of ordinary polytope $P=P^{d,k,n}$ such that the dimension of
$P$ is $d$, $P$ has $n+1$ vertices, and the first vertex of $P$ is on
exactly $k$ edges.
The vertex figure of the first vertex of $P^{d,k,n}$ is the cyclic 
$(d-1)$-polytope with $k$ vertices.
\end{theorem}

We use the following description of the facets of $P^{d,k,n}$ by 
Dinh \cite{Dinh}.
For any subset $X\subseteq {\mathbf Z}$, let $\mbox{ret}_n(X)$ 
(the ``retraction'' of $X$) be the set
obtained from $X$ by replacing every negative element by 0 and replacing every 
element greater than $n$ by $n$.
Let ${\cal X}_n$ be the collection of sets
\begin{equation} \label{facet}
X=[i,i+2r-1]\cup Y\cup [i+k,i+k+2r-1], \end{equation}
where $i\in{\bf Z}$, $1\le r\le m$, $Y$ is a paired $(d-2r-1)$-element
subset of $[i+2r+1,i+k-2]$, and $|\mbox{ret}_n(X)|\ge d$.
Dinh's theorem says that $\mbox{ret}_n({\cal X}_n)={\cal F}(P^{d,k,n})$,
the set of facets of $P^{d,k,n}$.
It is easy to check that when $n=k$, $|\mbox{ret}_n(X)|=d$ for all
$X\in{\cal X}_n$, and that $\mbox{ret}_n({\cal X}_n)$ is the set of
$d$-element Gale subsets of $[0,k]$, that is, the
facets of the cyclic polytope $P^{d,k,k}$.

Note that ${\cal X}_{n-1}\subseteq {\cal X}_n$.
We wish to describe ${\cal F}(P^{d,k,n})$ in terms of ${\cal F}(P^{d,k,n-1})$;
for this we need the following shift operations.
If $F=\mbox{ret}_{n-1}(X)\in{\cal F}(P^{d,k,n-1})$, let the right-shift of $F$
be $\mbox{rsh}(F)=\mbox{ret}_n(X+1)$.
Typically, $\mbox{rsh}(F)$ contains $F+1$.
If it does not, then $0\in F\setminus X$ and 
$\mbox{rsh}(F)=((F+1)\setminus\{1\})\cup \{0\}$.
In either case, $|\mbox{rsh}(F)|\ge |F|\ge d$.
If $F=\mbox{ret}_n(X)\in{\cal F}(P^{d,k,n})$, let the left-shift of $F$
be $\mbox{lsh}(F)=\mbox{ret}_{n-1}(X-1)$.
Note that $\mbox{lsh}(F)\setminus\{0\}=(F-1)\cap[1,n]$; $\mbox{lsh}(F)$ contains
0 if either 0 or 1 is in $F$.
\begin{lemma} If $n\ge k+1$ and $F\in{\cal F}(P^{d,k,n})$ with $\max F\ge k$,
then $\rm{lsh}(F)\in{\cal F}(P^{d,k,n-1})$.
\end{lemma}
\begin{proof}
Let $F=\mbox{ret}_n(X)$, with $X=[i,i+2r-1]\cup Y\cup [i+k,i+k+2r-1]$.
Then $X-1$ also has the form of equation~(\ref{facet}) (for $i-1$).
The set $\mbox{lsh}(F)$ is the vertex set of a facet of $P^{d,k,n-1}$ as long
as $|\mbox{lsh}(F)|\ge d$.
We check this in three cases.

\noindent{\em Case 1.}
If $k\le i+k+2r-1 \le n$, then $i+2r-1\ge 0$, so
$Y\subseteq[i+2r+1,i+k-2]\subseteq[2,i+k-2]$.
Then
$$\mbox{lsh}(F)\supseteq \max\{i+2r-2,0\}\cup 
(Y-1)\cup[i+k-1,i+k+2r-2],$$
so $|\mbox{lsh}(F)|\ge 1+(d-2r-1)+2r=d$.

\noindent{\em Case 2.}
If $i+k\ge n$, then $i\ge n-k\ge 1$.
Also, $|F|\ge d$ implies $\max Y\le n-1$.  So
$$\mbox{lsh}(F)=[i-1,i+2r-2]\cup(Y-1)\cup \{n-1\},$$ so 
$|\mbox{lsh}(F)|=2r+(d-2r-1)+1=d$.

\noindent{\em Case 3.}
If $i+k < n < i+k+2r-1$, then $i+2r-1\ge n-k \ge 1$, and 
$$F=[\max\{0,i\},i+2r-1]\cup Y
\cup [i+k,n],$$ so 
\begin{eqnarray*}
|F|&=&(i+2r-\max\{0,i\})+(d-2r-1)+(n-i-k+1)\\
&=&d+n-k-\max\{i,0\}\ge d+1.
\end{eqnarray*}
Then $|\mbox{lsh}(F)|\ge |F|-1\ge d$.

Thus, $\mbox{lsh}(F)$ is a facet of $P^{d,k,n-1}$.
\end{proof}

Identify each facet of the ordinary polytope $P^{d,k,n}$ with its ordered list
of vertices.
Then order the facets of $P^{d,k,n}$ in colex order.
This means, if $F=i_1i_2\ldots i_p$ and $G=j_1j_2\ldots j_q$, then 
$F\prec_c G$
if and only if for some $t\ge 0$, $i_{p-t}<j_{q-t}$ while for $0\le s<t$,
$i_{p-s}=j_{q-s}$.

\begin{lemma} If $n\ge k+1$ and $F_1$ and $F_2$ are facets of $P^{d,k,n}$ with 
$\max F_i\ge k$, then 
$F_1\prec_c F_2$ implies $\mbox{lsh}(F_1)\prec_c \mbox{lsh}(F_2)$.
\end{lemma}
\begin{proof}
Suppose $F_1\prec_c F_2$, and let $q$ be the maximum vertex in $F_2$ not
in $F_1$.
Then $\mbox{lsh}(F_1)\prec_c \mbox{lsh}(F_2)$ as long as $q\ge 2$,
for in that case $q-1\in \mbox{lsh}(F_2)\setminus \mbox{lsh}(F_1)$,
while 
$[q,n-1]\cap\mbox{lsh}(F_1)=[q,n-1]\cap\mbox{lsh}(F_2)$.
(If $q = 1$, then $q$ shifts to 0 in $\mbox{lsh}(F_2)$, but 0 may be
in $\mbox{lsh}(F_1)$ as a shift of a smaller element.)
So we prove $q\ge 2$.
Write $$F_2=\mbox{ret}_n([i,i+2r-1]\cup Y\cup [i+k,i+k+2r-1])$$ and
$$F_1=\mbox{ret}_n([i',i'+2r'-1]\cup Y'\cup [i'+k,i'+k+2r'-1]).$$

Since $\max F_2\ge k$, $i+2r-1\ge 0$, so $Y\cup[i+k,i+k+2r-1]\subseteq[2,n]$,  
Thus, if $q\in Y\cup[i+k,i+k+2r-1]$, then $q\ge 2$.
Otherwise $Y\cup [i+k,i+k+2r-1])= Y'\cup [i'+k,i'+k+2r'-1])$, but $Y\ne Y'$.
This can only happen when $Y\cup [i+k,i+k+2r-1])$ is an interval; in this case
$i+k+2r-1\ge n+1$.
Then $q=i+2r-1=(i+k+2r-1)-k\ge n+1-k\ge 2$.
\end{proof}

\begin{proposition}\label{recurs-facets}
Let $n\ge k+1$.  The facets of $P^{d,k,n}$ are
\begin{eqnarray*}
\lefteqn{\{F:\mbox{$F\in{\cal F}(P^{d,k,n-1})$ and $\max F\le n-2$}\}}\\
& &{}\cup \{\mbox{rsh}(F):\mbox{$F\in{\cal F}(P^{d,k,n-1})$ and $\max F\ge n-2$}\}.\end{eqnarray*}
\end{proposition}
\begin{proof}
If $\max X\le n-2$, then $\mbox{ret}_n(X)=\mbox{ret}_{n-1}(X)$;
in this case, letting $F=\mbox{ret}_n(X)$, $F\in{\cal F}(P^{d,k,n-1})$ if 
and only if $F\in{\cal F}(P^{d,k,n})$.
If $F\in{\cal F}(P^{d,k,n-1})$ with $\max F\ge n-2$, then 
$\mbox{rsh}(F)\in{\cal F}(P^{d,k,n})$ with $\max\mbox{rsh}(F)\ge n-1$.
Now suppose that $G=\mbox{ret}_n(X)\in{\cal F}(P^{d,k,n})$ with
$\max G\ge n-1$.
Let $F=\mbox{lsh}(G)=\mbox{ret}_{n-1}(X-1)\in{\cal F}(P^{d,k,n-1})$;
then $\max F \ge n-2$.
By definition, $\mbox{rsh}(F)=\mbox{ret}_n((X-1)+1)=\mbox{ret}_n(X)=G$.
\end{proof}

\begin{theorem}\label{colex-shell}
Let $F_1$, $F_2$, \ldots, $F_v$ be the facets of $P^{d,k,n}$ in colex order.
Then
\begin{enumerate}
\item $F_1$, $F_2$, \ldots, $F_v$ is a shelling of $P^{d,k,n}$.\label{is-shell}
\item For each $j$ there is a unique minimal face $G_j$
      of $F_j$ not contained in $\displaystyle\cup_{i=1}^{j-1}F_i$.\label{uniq-min}
\item For each $j$, $2\le j\le v-1$, $G_j$ contains the vertex of $F_j$ of
      maximum index, and is contained in the $d-1$ highest vertices of $F_j$.
      \label{G_j-right}
\item For each $j$, the face $G_j$ is a simplex and the interval $[G_j,F_j]$ is
a Boolean lattice.\label{Boolean}
\end{enumerate}
\end{theorem}
Note that this theorem is not saying that the faces of $P^{d,k,n}$ in the 
interval $[G_j,F_j]$ are all simplices.

\begin{proof}
We construct explicitly the faces $G_j$ in terms of $F_j$.
The reader may wish to refer to the example that follows the proof.

{\bf Cyclic polytopes.}
We start with the cyclic polytopes.
(For the cyclics, the theorem is generally known, or at least a shorter
proof based on \cite{billera-lee} is possible, but we will need the description 
of the faces $G_j$ later.)

Let $F_1$, $F_2$, \ldots, $F_v$ be the facets, in colex order, of 
$P^{d,k,k}$, the cyclic $d$-polytope with vertex set $[0,k]$.
Each facet $F_j$ can be written as $F_j=I_j^0\cup I_j^1\cup I_j^2\cup  \cdots 
\cup I_j^p\cup I_j^k$, where $I_j^0$ is the interval of $F_j$ containing 0,
if $0\in F_j$, and $I_j^0=\emptyset$ otherwise; 
$I_j^k$ is the interval of $F_j$ containing $k$,
if $k\in F_j$, and $I_j^k=\emptyset$ otherwise; 
and the $I_j^\ell$ are the other (even) intervals of $F_j$ with the
elements of $I_j^\ell$ preceding the elements of $I_j^{\ell+1}$.
(For example, in $P^{7,9,9}$, $F_6=\{0,1,2,4,5,7,8\}$,
$I_6^0=\{0,1,2\}$, 
$I_6^1=\{4,5\}$, 
$I_6^2=\{7,8\}$, and
$I_6^9=\emptyset$.)
For the interval $[a,b]$, write $E([a,b])$ for the integers in the even
positions in the interval, that is, $E([a,b])=[a,b]\cap 
\{a+2i+1: i\in{\bf N}\}$.
Let $G_j=\cup_{\ell=1}^p E(I_j^\ell)\cup I_j^k$.
Since $I_j^0=F_j$ if and only if $j=1$, $G_1=\emptyset$, and for all
$j>1$, $G_j$ contains the maximum vertex of $F_j$.
Since $F_j$ is a simplex, $[G_j,F_j]$ is a Boolean lattice.

To show  that 
$F_1$, $F_2$, \ldots, $F_v$ is a shelling of $P^{d,k,k}$ we show that
$G_j$ is not in a facet before $F_j$ and that every ridge of $P^{d,k,k}$ in 
$F_j$ that does not contain $G_j$ is contained in a previous facet.
For $j>0$ the face $G_j$ consists of the right end-set $I_j^k$ (if
 nonempty) and the set $\cup_{j=1}^p E(I_j^\ell)$ of singletons.
Note that $G_j$ is contained in the $d-1$ highest vertices of $F_j$
(equivalently, $G_j$ does not contain the minimum element of $F_j$)
unless $j=v$, in which case $G_v=F_v$.
Any facet $F$ of $P^{d,k,k}$ containing $G_j$ must satisfy Gale's
evenness condition and therefore must contain an integer adjacent to 
each element of $\cup_{j=1}^pE(I_j^\ell)$.
If any element of the form $\max I_j^\ell+1$ is in $F$, then $F$
occurs after $F_j$ in colex order.
This implies that any $F_i$ previous to $F_j$ and containing $G_j$
also contains $\cup_{\ell=1}^pI_j^\ell\cup I_j^k$.
But $F_j$ is the first facet in colex order that contains 
$\cup_{\ell=1}^pI_j^\ell\cup I_j^k$.
So $G_j$ is not in a facet before $F_j$.

Now let $g\in G_j$; we wish to show that $F_j\setminus \{g\}$ is in
a previous facet.
If $g\in E(I_j^\ell)$ for $\ell>0$, let 
$F=F_j\setminus\{g\}\cup \{\min I_j^\ell-1\}$.
Then $F$ satisfies Gale's evenness condition and is a facet before
$F_j$.
Otherwise $g\in I_j^k\setminus E(I_j^k)$; in this case let
$F=F_j\setminus\{g\}\cup \{\max I_j^0+1\}$ (where we let 
$\max I_j^0+1=0$ if $I_j^0=\emptyset$).
Again $F$ satisfies Gale's evenness condition and is a facet before
$F_j$.

Thus the colex order of facets is a shelling order for the cyclic 
polytope $P^{d,k,k}$, and we have an explicit description for the 
minimal new face $G_j$ as $F_j$ is shelled on.

{\bf General ordinary.}
Now we prove the theorem 
for general $P^{d,k,n}$ by induction on $n\ge k$, for fixed $k$.
Among the facets of $P^{d,k,n}$, first in colex order are those with
maximum vertex at most $n-2$.
These are also the first facets in colex order of $P^{d,k,n-1}$.
Thus the induction hypothesis gives us that this initial segment is a 
partial shelling of $P^{d,k,n}$, and that assertions 2--4 hold for
these facets.

{\bf Later facets.}
It remains to consider the facets of $P^{d,k,n}$ 
ending in $n-1$ or $n$.
These facets come from shifting facets of $P^{d,k,n-1}$ ending in $n-2$ or 
$n-1$.
Our strategy here will be to prove statement~\ref{uniq-min} of the theorem for 
these facets.
The intersection of $F_j$ with $\cup_{i=1}^{j-1}F_i$ is then the antistar
of $G_j$ in $F_j$, and so it is the union of $(d-2)$-faces that form an 
initial segment of a shelling of $F_j$.
This will prove that the colex order $F_1$, $F_2$, \ldots, $F_v$ is a shelling 
of $P^{d,k,n}$.

Note that there is nothing to show for the last facet of $P^{d,k,n}$
in colex order. 
It is $F_v=[n-d+1,n]$, and is the only facet (other than the first)
whose vertex set forms a single interval.
Assume from now on that $j$ is fixed, with $j\le v-1$.
Later we will describe recursively the minimal new face $G_j$ as $F_j$ is 
shelled on.
It will always be the case that $\max F_j\in G_j$.
We will prove that $G_j$ is truly a new face (is not contained in a previous
facet), and that every ridge not containing all of $G_j$ is contained in
a previous facet.

{\bf Ridges not containing the last vertex.}
It is convenient to start by showing that every ridge of $P^{d,k,n}$
contained in $F_j$ and not containing
$\max F_j$ is contained in an earlier facet.
This case does not use the recursion needed for the other parts of
the proof.
Write 
$$X=[i,i+2r-1]\cup Y\cup [i+k, i+k+2r-1]$$ and
$F_j=\mbox{ret}_n(X)=\{z_1,z_2,\ldots, z_p\}$
with $0\le z_1<z_2<\cdots<z_p\le n$.
The facet $F_j$ is a $(d-1)$-multiplex, so its
facets are of the form
$$F_j(\hat{z}_t)=\{z_\ell:\mbox{$1\le \ell\le p$, $0<|\ell-t|\le d-2$}\}$$
for $2\le t\le p-1$, $F_j(\hat{z}_1)=\{z_1,z_2,\ldots, z_{d-1}\}$, and
$F_j(\hat{z}_p)=\{z_{p-d+2},\ldots,z_{p-1}, z_p\}$.
If $F_j(\hat{z}_t)$ does not contain $\max F_j=z_p$, then $t\le p-d+1$ and
this implies $i\le z_t\le i+2r-1$.
Consider such a $z_t$.

{\bf The first ridge.}
For $t=1$, there are three cases to consider.

\noindent{\em Case 1.}
Suppose $z_1\ge 1$.
Then $F_j(\hat{z}_1)=[i,i+2r-1]\cup Y$.
Let $I$ be the right-most interval of $F_j(\hat{z}_1)$.
Let $Z=(I-k)\cup F_j(\hat{z}_1)$, and $F=\mbox{ret}_n(Z)$.
Since $i\ge 1$ and $\max F_j(\hat{z}_1)\le i+k-2$, the interval $I-k$ 
contributes at least one new element to $F$, so $|F|\ge d$.

\noindent{\em Case 2.}
Suppose $z_1=0$ and the right-most interval of $F_j(\hat{z}_1)$ is odd.
In this case the left-most interval of $F_j$ must also be odd, so $i<0$,
and $F_j(\hat{z}_1)$ contains $i+k$ but not $i+k-1$.
Let $F=F_j(\hat{z}_1)\cup\{i+k-1\}$.

\noindent{\em Case 3.}
Suppose $z_1=0$ and the right-most interval of $F_j(\hat{z}_1)$ is even
(and then so is the left-most interval).
Then $F_j(\hat{z}_1)=[0,i+2r-1]\cup Y\cup[i+k,k-1]$ (where the
last interval is empty if $i=0$).
Let $$F=F_j(\hat{z}_1)\cup\{i+2r\}=\{0\}\cup[1,i+2r]\cup Y\cup[i+k,k-1].$$
(When $i=0$ and $r=(d-1)/2$, this gives $F=[0,d-1]$.)

In all cases $F$ is a facet of $P^{d,k,n}$ containing $F_j(\hat{z}_1)$.
It does not contain $\max F_j$, so $F\prec_c F_j$.

{\bf Deleting a later vertex.}
Now assume $2\le t\le p-d+1$; then $z_t\ge \max\{i+1,1\}$.
Here $$F_j(\hat{z}_t)=[\max\{i,0\},z_t-1]\cup[z_{t+1},i+2r-1]\cup Y\cup[i+k,z_t-1+k],$$
and $|F_j(\hat{z}_t)|=z_t-\max\{i,0\}+d-2\ge d-1$.
Also note that 
$z_t-1+k$ is the $(d-2)$nd element of $\{z_1,z_2,\ldots, z_p\}$
after $z_t$, so $z_t-1+k=z_{t+d-2}<z_p=\max F_j$.

\noindent{\em Case 1.}
If $z_t-i$ is even, let $F=F_j(\hat{z}_t)\cup\{i+2r\}$.
Then $F=\mbox{ret}_n(Z)$, where
$$Z=[i,z_t-1]\cup [z_t+1,i+2r]\cup Y\cup [i+k, z_t-1+k],$$ and $|F|\ge d$.

\noindent{\em Case 2.}
If $z_t-i$ is odd and $\max ([i,i+2r-1]\cup Y) < i+k-2$, let $F=\mbox{ret}_n(Z)$, where
$$Z=[i-1,z_t-1]\cup [z_t+1,i+2r-1]\cup Y\cup [i+k-1,z_t-1+k].$$
Then $F\supseteq F_j(\hat{z}_t)\cup\{i+k-1\}$, so $|F|\ge d$.

\noindent{\em Case 3.}
Finally, suppose $z_t-i$ is odd and $\max Y = i+k-2$.
Let $[q,i+k-2]$ be the right-most interval of $Y$, and let
$F=\mbox{ret}_n(Z)$, where
$$Z=[q-k,z_t-1]\cup[z_t+1,i+2r-1]\cup(Y\setminus[q,i+k-2])\cup[q,z_t-1+k].$$
Then $F\supseteq F_j(\hat{z}_t)\cup\{i+k-1\}$, so $|F|\ge d$.

In all cases, $F$ is a facet of $P^{d,k,n}$ containing $F_j(\hat{z}_t)$ and 
$\max F_j \not\in F$, so $F$ occurs before $F_j$ in colex order.

{\bf Determining the minimal new face.}
We now describe the faces $G_j$ recursively.
(We are still assuming that $\max F_j\ge n-1$.)
Let $G$ be the face of $\mbox{lsh}(F_j)$ 
that is the minimal new face
when $\mbox{lsh}(F_j)$ is shelled on, in the colex shelling of the 
polytope $P^{d,k,n-1}$.
Let $G_j=G+1$; then $G_j$ contains $\max F_j$.
Also, $G_j$ is a subset of the last $d-1$ vertices of $F_j$, which form a facet 
of $F_j$; hence $G_j$ is a face of $F_j$.
For any facet $F_i$ of $P^{d,k,n}$, $G_j\subseteq F_i$ if and only if 
$G\subseteq \mbox{lsh}(F_i)$.
So by the induction hypothesis, $G_j$ is not contained in a facet
occurring before $F_j$ in colex order.

{\bf Ridges in previous facets.}
It remains to show that any ridge of 
$P^{d,k,n}$ contained in $F_j$ but not
containing all of $G_j$ is contained in a facet prior to $F_j$.
Note that we have already dealt with those ridges not containing $\max F_j$.
Now let $g\in G$, $g_j=g+1\in G_j$, and assume $g_j\ne \max F_j$.
The only ridge of $P^{d,k,n}$ contained in $F_j$, containing $\max F_j$,
and not containing $g_j$ is $F_j(\hat{g}_j)$.

Let $H$ be the unique ridge of $P^{d,k,n-1}$ in $\mbox{lsh}(F_j)$ containing 
$\max(\mbox{lsh}(F_j))$, but not containing $g$.
By the induction hypothesis, $H$ is contained in a facet $F$ of $P^{d,k,n-1}$
occurring before $\mbox{lsh}(F_j)$ in colex order.
Suppose $F_j(\hat{g}_j)$ is contained in a facet $F_\ell$ of $P^{d,k,n}$
occurring after $F_j$ in colex order.
Then $H$ is contained in $\mbox{lsh}(F_\ell)$.
Thus the ridge $H$ of $P^{d,k,n-1}$ is contained in three different facets:
$F$ (occurring before $\mbox{lsh}(F_j)$ in colex order), $\mbox{lsh}(F_j)$,
and $\mbox{lsh}(F_\ell)$
(occurring after $\mbox{lsh}(F_j)$ in colex order).
This contradiction shows that the ridge
$F_j(\hat{g}_j)$ can only be contained in a facet of $P^{d,k,n}$
occurring before $F_j$ in colex order.

{\bf Boolean intervals.}
Finally to verify assertion 4 of 
the theorem, observe that every facet
$F_j$ is a $(d-1)$-dimensional multiplex.
The vertices of $G_j$ are among the last $d$ vertices of $F_j$ and
so are affinely independent \cite{Bisz-mult}; thus $G_j$ is a simplex.
The maximum vertex $u$ of $F_j$ is contained in $G_j$.
The vertex figure of the maximum vertex in any multiplex is a simplex
\cite{Bisz-mult}.
So the interval $[G_j,F_j]$, which is an interval in $[u,F_j]$, 
is a Boolean lattice.
\end{proof}

A nonrecursive description of the faces $G_j$, generalizing that for the
cyclic case in the proof, is as follows.
Write the facet $F_j$ as a disjoint union, 
$F_j=A_j^0\cup I_j^1\cup I_j^2\cup  \cdots \cup I_j^p\cup I_j^n$, where 
$I_j^n$ is the interval of $F_j$ containing $n$ if $n\in F_j$, and
$I_j^n=\emptyset$ otherwise;
the $I_j^\ell$ ($1\le \ell\le p$) are even intervals of $F_j$ written in 
increasing order; and
$A_j^0$ is 
\begin{itemize}
\item the interval containing 0, if $\max F_j\le k-1$;
\item the union of the interval containing $\max F_j-k$ and the interval
      containing $\max F_j-k+2$ (if the latter exists), if 
      $k\le \max F_j\le n-1$;
\item the interval containing $n-k$, if $\max F_j=n$ and $n-k\in F_j$;
\item $\emptyset$, if $\max F_j=n$ and $n-k\not\in F_j$.
\end{itemize}
Then $G_j=\cup_{\ell=1}^p E(I_j^\ell)\cup I_j^n$.

\vspace{\baselineskip}

{\bf Example.}
Table \ref{P568} gives the faces $F_j$ and $G_j$ for the colex shelling of the
ordinary polytope $P^{5,6,8}$.

\begin{table}[hbt]
\begin{center}
\begin{tabular}{|l|l|r||l|l|r|}
\hline
$j$ & \multicolumn{1}{c|}{$F_j$} & 
\multicolumn{1}{c||}{$G_j$} & $j$ & 
\multicolumn{1}{c|}{$F_j$} & \multicolumn{1}{c|}{$G_j$}\\
\hline
1 & \verb[01234    [ & {\tt $\emptyset$} &
9 & \verb[  23 56 8[ & \verb[68[\\
2 & \verb[012 45   [ & \verb[5[ &
{10} & \verb[   3456 8[ & \verb[468[\\
3 & \verb[0 2345   [ & \verb[35[ &
{11} & \verb[ 1234  78[ & \verb[78[ \\
4 & \verb[0 23 56  [ & \verb[6[ &
{12} & \verb[ 12 45 78[ & \verb[578[\\
5 & \verb[0  3456  [ & \verb[46[ & 
{13} & \verb[0123  678[ & \verb[678[\\
6 & \verb[01 34 67 [ & \verb[7[ &
{14} & \verb[   34 678[ & \verb[4678[\\
7 & \verb[01  4567 [ & \verb[57[ &
{15} & \verb[012  5678[ & \verb[5678[\\
8 & \verb[  2345  8[ & \verb[8[ &
{16} & \verb[    45678[ & \verb[45678[\\
\hline
\end{tabular}
\end{center}
\caption{Shelling of $P^{5,6,8}$}\label{P568}
\end{table}

Let us look at what happens when facet $F_{13}$ is shelled on.
The ridges of $P^{5,6,8}$ contained in $F_{13}$ are 0123, 0236,
01367, 012678, 12378, 2368, and 3678.
The first ridge, 0123, is contained in $F_1=01234$.
The ridge 0236 is $F_{13}(\hat{z}_2)=F_{13}(\hat{1})$, and
$\max([i,i+2r-1]\cup Y)=3 < 4=i+k-2$, so we find that 0236 is
contained in $F_4 = 02356$.
The ridge 01367 is $F_{13}(\hat{z}_3)=F_{13}(\hat{2})$, so we
find that 01367 is contained in $F_6=013467$.
This facet $F_{13}=0123678$ is shifted from the facet 012567 of
$P^{5,6,7}$, which in turn is shifted from the facet 01456 of the
cyclic polytope $P^{5,6,6}$.
When 01456 occurs in the shelling of the cyclic polytope, its 
minimal new face is its right interval, 456.
In $P^{5,6,8}$, then, 
the minimal new face when $F_{13}$ is shelled on is 678.
The other ridges of $F_{13}$ not containing 678 are 12378 and
2368.
The interval $[G_{13},F_{13}]$ contains the triangle
678, the 3-simplex 3678, the 3-multiplex 012678, and $F_{13}$
itself (which is a pyramid over 012678).

\vspace{\baselineskip}

Note that for the multiplex, $M^{d,n}=P^{d,d,n}$, this theorem gives
a shelling different from the one mentioned in Section~\ref{triang-multi}.
In the standard notation for the facets of the multiplex 
(see Definition~\ref{def-multi}), 
the colex shelling order is $F_0$, $F_1$, \ldots, $F_{n-d}$, $F_{n-1}$,
$F_{n-2}$, \ldots, $F_{n-d+1}$, $F_n$.
The statements of this section hold also for even-dimensional multiplexes.

\section{The $h$-vector from the shelling}
The $h$-vector of a simplicial polytope can be obtained easily from
any shelling of the polytope.
For $P$ a simplicial polytope, and $\cup [G_j,F_j]$ the partition of a
face lattice of $P$ arising from a shelling, $h(P,x)=\sum_j x^{d-|G_j|}$.
For general polytopes, the (toric) $h$-vector can also be decomposed according 
to the shelling partition.
For a shelling, $F_1$, $F_2$, \ldots, $F_n$, of a polytope $P$, write
${\cal G}_j$ for the set of faces of $F_j$ not in $\cup_{i<j} F_i$.
Then $h(P,x)=\sum_{j=1}^n h({\cal G}_j,x)$, where
$h({\cal G}_j,x)=\sum_{G\in{\cal G}_j} g(G,x)(x-1)^{d-1-\dim G}$.
However, in general we do not know that the coefficients of $h({\cal G}_j,x)$
count anything natural, nor even that they are nonnegative.
Stanley raised this issue in \cite[Section 6]{sta87}.  It has apparently
been settled by Tom Braden \cite{braden}.

We turn now to $h$-vectors of ordinary polytopes.
In \cite{bayer-ordflags} we used the flag vector of the
ordinary polytope to compute its toric $h$-vector.
\begin{theorem}[\cite{bayer-ordflags}] For $ n\ge k\ge d=2m+1\ge 5$ and 
$1\le i\le m$,
$$h_i(P^{d,k,n})={k-d+i\choose i}+(n-k){k-d+i-1\choose i-1}.$$
\end{theorem}
We did not understand why the $h$-vector turned out to have such a nice
form.
Here we show how the $h$-vector can be computed from the colex shelling.
Properties~2 and~4 of Theorem~\ref{colex-shell} are critical.

In \cite{bayer-ordflags} we showed that the flag vector of a multiplicial
polytope depends only on the $f$-vector.
However, for our purposes here it is more useful to write 
the $h$-vector in terms of the $f$-vector and the flag vector entries of
the form $f_{0i}$ (the sum of the number of vertices on the $i$-faces)
We introduce a modified $f$-vector.
Let $\bar{f}_{-1}=f_{-1}=1$, $\bar{f}_0=f_0$, and
$\bar{f}_{d-1}=f_{d-1}+(f_{0,d-1}-df_{d-1})$; and for $1\le j\le d-2$, let
$$\bar{f}_j=f_j+(f_{0,j+1}-(j+2)f_{j+1})+(f_{0,j}-(j+1)f_j).$$
\begin{theorem}
If $P$ is a multiplicial $d$-polytope, then
$$h(P,x)=\sum_{i=0}^d h_i(P) x^{d-i}=\sum_{i=0}^d \bar{f}_{i-1}(P)(x-1)^{d-i}.$$
\end{theorem}
\begin{proof}
As observed in the proof of Theorem~\ref{trian-multi}, the $g$-polynomial of an 
$e$-dimensional multiplex $M$ with $n+1$ vertices is $g(M,x)=1+(n-e)x$.
So for a multiplicial $d$-polytope $P$,
\begin{eqnarray*}
h(P,x)&=& \sum_{{\mbox{{\scriptsize $G$ face of $P$}}}\atop
       {G\not=P}} g(G,x)(x-1)^{d-1-\dim G}\\
      &=& \sum_{{\mbox{{\scriptsize $G$ face of $P$}}}\atop
       {G\not=P}} (1+(f_0(G)-1-\dim G)x)(x-1)^{d-1-\dim G}\\
      &=&\sum_{i=0}^d f_{i-1}(x-1)^{d-i}
        +\sum_{i=1}^{d-1}(f_{0i}-(i+1)f_i)x(x-1)^{d-1-i}\\
      &=&\sum_{i=0}^d f_{i-1}(x-1)^{d-i}
        +\sum_{i=1}^{d-1}(f_{0i}-(i+1)f_i)[(x-1)^{d-i}+(x-1)^{d-1-i}]\\
      &=&(x-1)^d+f_0(x-1)^{d-1}\\
      & &{}+\sum_{i=2}^{d-1}(f_{i-1}+(f_{0i}-(i+1)f_i)
        +(f_{0,i-1}-if_{i-1}))(x-1)^{d-i}\\
      & &{}+(f_{d-1}+(f_{0,d-1}-df_{d-1}))\\
      &=&\sum_{i=0}^d \bar{f}_{i-1}(P)(x-1)^{d-i}.
\end{eqnarray*}
\end{proof}

Simplicial polytopes are a special case of multiplicial polytopes.
Clearly, when $P$ is simplicial, $\bar{f}(P)=f(P)$, and we recover the
definition of the simplicial $h$-vector in terms of the $f$-vector.
The multiplicial $h$-vector formula can be thought of as breaking
into two parts: one involving the $f$-vector, and matching the
simplicial $h$-vector formula; the other involving the ``excess
vertex counts,'' $f_{0,j}-(j+1)f_j$.
In the simplicial case the sum of the entries in the $h$-vector is
the number of facets.
For multiplicial polytopes 
$\sum_{i=0}^d h_i(P)=\bar{f}_{d-1}(P)=f_{d-1}+
(f_{0,d-1}-df_{d-1})$.

In general, applying the simplicial $h$-formula to a nonsimplicial
$f$-vector produces a vector with no (known) combinatorial interpretation.
This vector is neither symmetric nor nonnegative in general.
We will see that in the case of ordinary polytopes something special happens.
Write $h'(P,x)=\sum_{i=0}^d h'_i(P)x^{d-i}$ for the $h$-polynomial that $P$ 
would have if it were simplicial,
by letting $h'_i(P)$ be the coefficient of $x^{d-i}$ in 
$\sum_{i=0}^d f_{i-1}(P)(x-1)^{d-i}$.

\begin{theorem}\label{h'}
Let $P^{d,k,n}$ be an ordinary polytope. 
Let $\cup_{j=1}^v[G_j,F_j]$ be the partition of the face lattice of
$P^{d,k,n}$ associated with the colex shelling of $P^{d,k,n}$.
Then for all $i$, $0\le i\le d$, 
$h'(P^{d,k,n},x)= \sum_{j=1}^v x^{d-|G_j|}$.

Furthermore, if $C^{d,k}$ is the cyclic $d$-polytope with $k+1$ vertices, 
then for all $i$, $0\le i\le d$, $h'_i(P^{d,k,n})\ge h_i(C^{d,k})$, with 
equality for $i > d/2$.
\end{theorem}
\begin{proof}
Direct evaluation gives $h'_0(P)=h'_d(P)=1$.
Let $F_1$, $F_2$, \ldots, $F_v$ be the colex shelling of $P^{d,k,n}$.
By Theorem~\ref{colex-shell}, part~2, the set of faces of $P^{d,k,n}$
has a partition as $\cup_{j=1}^v [G_j,F_j]$.
By Theorem~\ref{colex-shell}, part~4, the interval $[G_j,F_j]$ has exactly
${d-1-\dim G_j \choose \ell-\dim G_j}$ faces of dimension $\ell$\/ for
$\dim G_j\le \ell \le d-1$.
Let $k_i=|\{j:\dim G_j=i-1\}|=|\{j: |G_j|=i\}|$ (since all the $G_j$ are
simplices).
Then $f_\ell=\sum_{i=0}^{l+1}{d-i \choose \ell-i+1} k_i$.
These are the (invertible) equations that give $f_\ell$ in terms of $h'_i$,
so for all $i$, $h'_i=k_i= |\{j: |G_j|=i\}|$.

The second part we prove by induction on $n\ge k$.
Recall from the proof of Theorem~\ref{colex-shell} that for each facet $F_j$ 
of $P^{d,k,n}$, $G_j$ is the same size as the minimum new face $G$ of 
the corresponding facet of $P^{d,k,n-1}$; that facet is the same (as vertex
set) as $F_j$, if
$\max F_j\le n-2$, and is $\mbox{lsh}(F_j)$, if $\max F_j\ge n-1$.
From Proposition~\ref{recurs-facets} we see that each facet of $P^{d,k,n-1}$
with maximum vertex $n-2$ gives rise to two facets of $P^{d,k,n}$,
while all others give rise to exactly one facet each.
So for all $i$, 
\begin{eqnarray*}
\lefteqn{h'_i(P^{d,k,n})=h'_i(P^{d,k,n-1})}\\
&+&
|\{j:\mbox{$F_j$ is a facet of $P^{d,k,n}$ with $\max F_j=n-1$ and 
$|G_j|=i$}\}|.
\end{eqnarray*}
Thus, $h'_i(P^{d,k,n})\ge h'_i(P^{d,k,n-1})$, so by induction,
$h'_i(P^{d,k,n})\ge h'_i(C^{d,k})$ for all $i$.  

Now consider $F_j$ with $\max F_j=n-1$.
In the description of the faces $G_j$ following the proof of 
Theorem~\ref{colex-shell}, $I_j^n=\emptyset$ and 
$F_j\setminus A_j^0\subseteq Y\cup [i+k,i+k+2r-1]$.
Thus $|G_j|\le |E(Y\cup [i+k,i+k+2r-1])|=(d-1)/2$.
So for $i>d/2$, $h'_i(P^{d,k,n})=h'_i(P^{d,k,n-1})=h_i(C^{d,k})$.
\end{proof}

Note that for the multiplex $M^{d,n}$ ($d$ odd or even), 
$h'(M^{d,n})=(1,n-d+1,1,1,\ldots, 1, 1)$, while
$h(M^{d,n})=(1,n-d+1,n-d+1,\ldots, n-d+1, 1)$.

Now for multiplicial polytopes, we consider the remaining part of the 
$h$-vector, coming from the parameters $f_{0,j}-(j+1)f_j$.
This is 
\begin{eqnarray*}
\lefteqn{h(P,x)-h'(P,x)}\\
&=&(f_{0,d-1}-df_{d-1})+\sum_{i=2}^{d-1}\left(
(f_{0,i}-(i+1)f_i)+(f_{0,i-1}-if_{i-1})\right)(x-1)^{d-i}.
\end{eqnarray*}
So
\begin{eqnarray*}
\lefteqn{h(P,x+1)-h'(P,x+1)}\\
&=&(f_{0,d-1}-df_{d-1})+\sum_{i=2}^{d-1}
\left((f_{0,i}-(i+1)f_i)+(f_{0,i-1}-if_{i-1})\right)x^{d-i}\\
&=&\sum_{i=2}^{d-1} (f_{0,i}-(i+1)f_i)(x+1)x^{d-1-i}.
\end{eqnarray*}
So $$ \sum_{i=2}^{d-1} (h_i(P)-h'_i(P))(x+1)^{d-1-i}=
\sum_{i=2}^{d-1} (f_{0,i}-(i+1)f_i)x^{d-1-i}.$$

For the ordinary polytope,
this equation can be applied locally to give the contribution
to $h(P^{d,k,n},x)-h'(P^{d,k,n},x)$ from each interval $[G_j,F_j]$ of the 
shelling partition.
For each $j$, and each $i\ge \dim G_j$, let 
$b_{j,i}=\sum(f_0(H)-(i+1))$, where the sum is over all $i$-faces $H$
in $[G_j,F_j]$.
Let $b_j(x)=\sum_{i=\dim G_j}^{d-1} b_{j,i}x^{d-1-i}$.
Write $b_j(x)$ in the basis of powers of $(x+1)$:
$b_j(x)=\sum a_{j,i}(x+1)^{d-1-i}$.
Then $a_{j,i}=h_i({\cal G}_j)-h'_i({\cal G}_j)$, the contribution to 
$h_i(P^{d,k,n})-h'_i(P^{d,k,n})$ from faces in the interval $[G_j,F_j]$.
Note that for fixed $j$, $\sum_i a_{j,i}=b_j(0) = f_0(F_j)-d$.
We will return to the coefficients $a_{j,i}$ after triangulating
the ordinary polytope.

\vspace{\baselineskip}

{\bf Example.}
The $h$-vector of $P^{5,6,8}$ is $h(P^{5,6,8})=(1,4,7,7,4,1)$.
The sum of the $h_i$ is 24, which counts the 16 facets plus one
for each of the four 6-vertex facets, plus two for each of the
two 7-vertex facets.
Referring to Table~1, we see that
$h'(P^{5,6,8})=(1,4,5,3,2,1)$; from this we compute 
$f(P^{5,6,8})=(9,31, 52, 44, 16)$.
The nonzero $a_{j,i}$ here are $a_{6,2}=a_{7,3}=a_{11,2}=a_{12,3}=1$
and $a_{13,3}=a_{15,4}=2$.
In this case each interval $[G_j,F_j]$ contributes to 
$h_i(P^{d,k,n})-h'_i(P^{d,k,n})$ for at most one $i$, but this is
not true in general.

\section{Triangulating the ordinary polytope}
Triangulations of polytopes or of their boundaries can be used
to calculate the $h$-vector of the polytope if the triangulation
is shallow \cite{Bayer-weakly}.
The solid ordinary polytope need not have a shallow triangulation,
but its boundary does have a shallow triangulation.
The triangulation is obtained simply by triangulating each multiplex
as in Section~\ref{triang-multi}.
This triangulation is obtained by pushing the vertices in the order
0, 1, \ldots, $n$.
\begin{theorem}
The boundary of the ordinary polytope $P^{d,k,n}$ has a shallow triangulation.
The facets of one such triangulation are the Gale subsets of $[i,i+k]$
(where $0\le i\le n-k$)
of size $d$ containing either 0 or $n$ or the set $\{i,i+k\}$.
\end{theorem}
\begin{proof}
First we show that each such set is a consecutive subset of some facet of 
$P^{d,k,n}$.
Suppose $Z$ is a Gale subset of $[i,i+k]$ of size $d$ containing $\{i,i+k\}$.
Write $Z=[i,i+a-1]\cup Y\cup [i+k-b+1,i+k]$,
where $a\ge 1$, $b\ge 1$, and $Y\cap \{i+a,i+k-b\}=\emptyset$.
Since $Z$ is a Gale subset, $|Y|$ is even; let $r=(d-1-|Y|)/2$.
Since $|Z|=d$, $a+b=2r+1$, so $a$ and $b$ are each at most $2r$.
Define $X=[i+a-2r,i+a-1]\cup Y\cup [i+k-b+1,i+k-b+2r]$.  
Note that $i+k-b+1=(i+a-2r)+k$.
Then $\mbox{ret}_n(X)$ is the vertex set of a facet of $P^{d,k,n}$, and 
$Z$ is a consecutive subset of $\mbox{ret}_n(X)$.

Now suppose that $Z$ is a Gale subset of $[0,k]$ of size $d$ containing 0,
but not $k$.
Write $Z=\{0\}\cup Y\cup [j-2r+1,j]$, where $j<k$, $r\ge 1$, 
and $j-2r\not\in Y$.
Then $|Y|=d-2r-1$, and $Z=\mbox{ret}_n(X)$, where 
$X=[j-2r+1-k,j-k]\cup Y\cup [j-2r+1,j]$.
So $Z$ itself is the vertex set of a facet of $P^{d,k,n}$.
The case of sets containing $n$ but not $n-k$ works the same way.

Next we show that all consecutive $d$-subsets of facets $F$ of $P^{d,k,n}$ are 
of one of these types.
Let $F=\mbox{ret}_n(X)$, where $X=[i,i+2r-1]\cup Y\cup [i+k,i+k+2r-1]$,
with $Y$ a paired subset of size $d-2r-1$ of $[i+2r+1,i+k-2]$.
Suppose first that $i+2r-1\ge 0$ and $i+k\le n$.
Let $Z$ be a consecutive $d$-subset of $F$.
Since $|Y|=d-2r-1$, $|[i,i+2r-1]\cap F|\le 2r$, and 
$|[i+k,i+k+2r-1]\cap F|\le 2r$, it follows that 
$i+2r-1$ and $i+k$ must both be in $Z$.
Thus we can write $Z=[i+2r-a,i+2r-1]\cup Y \cup
[i+k,i+k+b-1]$, with $a+b=2r+1$, $i+2r-a\ge 0$, 
and $i+k+b-1\le n$.
Let $\ell=i+2r-a$. 
Then $i+k+b-1=\ell+k$, so $0\le \ell\le n-k$, and 
$Z$ is a Gale subset of $[\ell,\ell+k]$ containing $\{\ell,\ell+k\}$.

If $i+2r-1<0$, then $i+k+2r-1<k\le n$, and $F=\{0\}\cup Y\cup[i+k,i+k+2r-1]$.
Then $|F|=d$ and $F$ itself is a Gale subset of $[0,k]$ of size $d$
containing 0.
Similarly for the case $i+k > n$.

The sets described are exactly the $(d-1)$-simplices obtained by 
triangulating each facet of $P^{d,k,n}$ according to 
Theorem~\ref{trian-multi}.
The fact that this triangulation is shallow follows from the corresponding
fact for the triangulations of the multiplexes.
\end{proof}

Let ${\cal T}={\cal T}(P^{d,k,n})$ be this triangulation of 
$\partial P^{d,k,n}$.
Since $\cal T$ is shallow, $h(P^{d,k,n},x)=h({\cal T},x)$.
We calculate $h({\cal T},x)$ by shelling $\cal T$.

\begin{theorem}\label{triang-shell}
Let $F_1$, $F_2$, \ldots, $F_v$ be the colex order of the facets of $P^{d,k,n}$.
For each $j$, if $F_j=\{z_1,z_2,\ldots, z_{p_j}\}$ ($z_1<z_2<\cdots<z_{p_j}$),
and $1\le \ell\le p_j-d+1$, let 
$T_{j,\ell}=\{z_\ell,z_{\ell+1}, \ldots, z_{\ell+d-1}\}$.
Then $T_{1,1}$, $T_{1,2}$, \ldots, $T_{1,p_1-d+1}$, $T_{2,1}$, \ldots, 
$T_{2,p_2-d+1}$, \ldots, $T_{v,1}$, \ldots, $T_{v,p_v-d+1}$ is a shelling of 
${\cal T}(P^{d,k,n})$.

Let $U_{j,\ell}$ be the minimal new face when $T_{j,\ell}$ is shelled on.
As vertex sets, $U_{j,p_j-d+1}=G_j$.
\end{theorem}
\begin{proof}
Throughout the proof, write $F_j= \{z_1,z_2,\ldots, z_{p_j}\}$ 
($z_1<z_2<\cdots<z_{p_j}$).
We first show that $G_j$ is the unique minimal face of $T_{j,p_j-d+1}$ not
contained in 
$\displaystyle(\cup_{i=1}^{j-1}\cup_{\ell=1}^{p_i-d+1}T_{i,\ell})\cup
(\cup_{\ell=1}^{p_j-d}T_{j,\ell})$.
The set $G_j$ is not contained in a facet of $P^{d,k,n}$ earlier than $F_j$.
So $G_j$ does not occur in a facet of $\cal T$ of the form $T_{i,\ell}$ for
$i<j$.
Also, $\max F_j\in G_j$, so $G_j$ does not occur in a facet of $\cal T$ of the
form $T_{j,\ell}$ for $\ell\le p_j-d$.
Thus $G_j$ does not occur in a facet of $\cal T$ before $T_{j,p_j-d+1}$.

We show that for $z_q\in G_j$, $T_{j,p_j-d+1}\setminus \{z_q\}$ is contained
in a facet of $\cal T$ occurring before $T_{j,p_j-d+1}$.
There is nothing to check for $j=v$, because $T_{v,p_v-d+1}=F_v$
is the last simplex in the purported shelling order.
So we may assume that $j<v$ and thus $G_j$ is contained in the last $d-1$
vertices of $F_j$.

\noindent{\em Case 1.} If $p_j>d$ and $q=p_j$ ($z_q$ is the maximal element of
$F_j$), then 
$T_{j,p_j-d+1}\setminus \{z_{p_j}\}\subset T_{j,p_j-d}$.

\noindent{\em Case 2.} Suppose $p_j-d+2\le q\le p_j-1$.
Then $T_{j,p_j-d+1}\setminus \{z_q\}\subseteq \{z_{q-d+2},\ldots, z_{q-1},
z_{q+1}\ldots, z_{p_j}\}=H$. 
This is a ridge of $P^{d,k,n}$ in $F_j$
not containing $G_j$,
and hence $H$ is contained in a previous facet $F_\ell$ of $P^{d,k,n}$.
Since $H$ is a ridge in both $F_j$ and $F_\ell$, $H$ is obtained from each
facet by deleting a single element from a consecutive string of vertices in
the facet.
So $|H|\le |F_\ell\cap[z_{q-d+2},z_{p_j}]|\le |H|+1$, and so
$d-1\le |F_\ell\cap[z_{p_j-d+1},z_{p_j}]|\le d$.
So $T_{j,p_j-d+1}\setminus \{z_q\}$ is contained in 
a consecutive set of
$d$ elements of $F_\ell$, and hence in a $(d-1)$-simplex of 
${\cal T}(P^{d,k,n})$ belonging to $F_\ell$.
This simplex occurs before $T_{j,p_j-d+1}$ in the specified shelling order.

\noindent{\em Case 3.} Otherwise $p_j=d$ (so $p_j-d+1=1$) and $q=d$. 
Then $T_{j,1}=F_j$
and $H=T_{j,1}\setminus\{z_d\}$ is a ridge of $P^{d,k,n}$ in $F_j$ not
containing $\max F_j$, so $H$ is contained in a previous facet 
$F_\ell$ of $P^{d,k,n}$.
As in Case 2, 
$d-1\le |F_\ell\cap[z_1,z_{d-1}]|\le d$.
So $T_{j,1}\setminus\{z_d\}$ is contained in a consecutive set of
$d$ elements of $F_\ell$, and hence in a $(d-1)$-simplex of 
${\cal T}(P^{d,k,n})$ belonging to $F_\ell$.
This simplex occurs before $T_{j,p_j-d+1}$ in the specified shelling order.

So in the potential shelling of $\cal T$, $G_j$ is the unique minimal
new face as $T_{j,p_j-d+1}$ is shelled on.
Write $U_{j,p_j-d+1}=G_j$.
At this point we need a clearer view of the simplex $T_{j,\ell}$.
Recall that $F_j$ is of the form $\mbox{ret}_n(X)$, where 
$X=[i,i+2r-1]\cup Y\cup [i+k,i+k+2r-1]$, with $Y$ a subset of size $d-2r-1$.
If $i+2r-1<0$ or $i+k>n$, then $p_j=|F_j|=d$, and $T_{j,1}=T_{j,p_j-d+1}=F_j$;
we have already completed this case.
So assume $i+2r-1\ge 0$ and $i+k\le n$.
A consecutive string of length $d$ in $\mbox{ret}_n(X)$ must then be of the
form $[i+s,i+2r-1]\cup Y\cup[i+k,i+k+s]$ for some $s$, $0\le s\le 2r-1$.
(All such strings---with appropriate $Y$---having $i+s\ge 0$ and 
$i+k+s\le n$ occur as $T_{j,\ell}$.)
In particular, for 
$\ell<p_j-d+1$, $T_{j,\ell}=T_{j,\ell+1}\setminus\{\max T_{j,\ell+1}\}\cup
\{\min T_{j,\ell+1}-1\}$ and
$\max T_{j,\ell}=\min T_{j,\ell}+k$.

Now define $U_{j,\ell}$ for $\ell\le p_j-d$ recursively by 
$U_{j,\ell}=U_{j,\ell+1}\setminus\{z\}\cup \{z-k,z-1\}$, where
$z=\max T_{j,\ell+1}$.
By the observations above, $U_{j,\ell}\subseteq T_{j,\ell}$.
We prove by downward induction that $U_{j,\ell}$ is not contained in a facet
$F_i$ of $P^{d,k,n}$ before $F_j$, that
$U_{j,\ell}$ is not contained in a facet
of $\cal T$ occurring before $T_{j,\ell}$, and that any ridge of $\cal T$
in $T_{j,\ell}$ not containing all of $U_{j,\ell}$ is in an earlier facet of
$\cal T$.
The base case of the induction is $\ell=p_j-d+1$, and this case has been
handled above.

Note that $\{z-k,z-1\}$ is a diagonal of the 2-face
$\{z-k-1,z-k,z-1,z\}$ of $P^{d,k,n}$ \cite{Bisz}.
So if $F_i$ is a facet of $P^{d,k,n}$ containing $U_{j,\ell}$, then $F_i$
contains $\{z-k-1,z-k,z-1,z\}$.
Thus $F_i$ contains $U_{j,\ell+1}$, so, by the induction assumption, $i\ge j$.
Therefore, for $i<j$, and any $r$, $T_{i,r}$ does not contain $U_{j,\ell}$.
For $r<\ell$, $T_{j,r}$ does not contain $z-1=\max T_{j,\ell}$, 
so $T_{j,r}$ does not contain $U_{j,\ell}$.

Now we wish to show that for any $g\in U_{j,\ell}$, $T_{j,\ell}\setminus \{g\}$ 
is in a previous facet of $\cal T$.

\noindent {\em Case 1.}
If $g=z-1=\max T_{j,\ell}$ and $\ell\ge 2$, then $T_{j,\ell}\setminus\{g\}
\subset T_{j,\ell-1}$.

\noindent{\em Case 2.}
If $g=z-1=\max T_{j,\ell}$ and $\ell = 1$, then $T_{j,\ell}\setminus\{g\}$
is the leftmost ridge of $P^{d,k,n}$ in $F_j$ and, in particular, does not
contain $\max F_j$.
So $H=T_{j,\ell}\setminus\{g\}$ is contained in a previous facet $F_e$
of $P^{d,k,n}$.
As in the $\ell=p_j-d+1$ case, $F_e\cap [\min T_{j,\ell}, \max T_{j,\ell}]$
is contained in a consecutive set of $d$ elements of $F_e$, and hence in a
$(d-1)$-simplex of ${\cal T}(P^{d,k,n})$ belonging to $F_e$.
So $T_{j,\ell}\setminus\{g\}$ is contained in a previous facet of $\cal T$.

\noindent{\em Case 3.}
Suppose $g<z-1$ and $g\in U_{j,\ell}\cap U_{j,\ell+1}$.
Since $\{z-1,z\}\subset T_{j,\ell+1}$, $T_{j,\ell+1}$ contains at most
$d-3$ elements less than $g$.
The ridge $H$ of $P^{d,k,n}$ in $F_j$ containing $T_{j,\ell+1}\setminus\{g\}$
consists of the $d-2$ elements of $F_j$ below $g$ and the (up to) $d-2$ 
elements of $F_j$ above $g$.
In particular, $H$ contains $\min T_{j,\ell+1}-1=\min T_{j,\ell}$.
So $T_{j,\ell}\setminus\{g\}\subset H$.
Since $\dim T_{j,\ell}\setminus\{g\}=d-2$, $H$ is the (unique) smallest face of 
$P^{d,k,n}$ containing $T_{j,\ell+1}\setminus\{g\}$.
By the induction hypothesis $T_{j,\ell+1}\setminus\{g\}$ is contained in
a previous facet $T_{i,r}$ of $\cal T$; here $i<j$ because
$\max T_{j,\ell+1}\in T_{j,\ell+1}\setminus\{g\}$.
The $(d-2)$-simplex $T_{j,\ell+1}\setminus\{g\}$ is then contained in a
ridge of $P^{d,k,n}$ contained in $F_i$, but this ridge must be $H$,
by the uniqueness of $H$.
So $T_{j,\ell}\setminus\{g\}\subset H=F_i\cap F_j$.
As in earlier cases, $F_i\cap [\min T_{j,\ell}, \max T_{j,\ell}]$
is contained in a consecutive set of $d$ elements of $F_i$, and hence in a
$(d-1)$-simplex of ${\cal T}(P^{d,k,n})$ belonging to $F_i$.
So $T_{j,\ell}\setminus\{g\}$ is contained in a previous facet of $\cal T$.

\noindent{\em Case 4.}
Finally, let $g=z-k$, which is $\min T_{j,\ell}+1$.
Then $T_{j,\ell}$ contains $d-2$ elements above $g$.
Let $H$ be the ridge of $P^{d,k,n}$ in $F_j$ containing 
$T_{j,\ell}\setminus \{g\}$.
Then $\max H=\max T_{j,\ell}<\max F_j$, so $H$ does not contain $G_j$.
So $H$ is in a previous facet $F_i$ of $P^{d,k,n}$.
As in earlier cases, $F_i\cap [\min T_{j,\ell}, \max T_{j,\ell}]$
is contained in a consecutive set of $d$ elements of $F_i$, and hence in a
$(d-1)$-simplex of ${\cal T}(P^{d,k,n})$ belonging to $F_i$.
So $T_{j,\ell}\setminus\{g\}$ is contained in a previous facet of $\cal T$.

Thus $T_{1,1}$, $T_{1,2}$, \ldots, $T_{1,p_1-d+1}$, $T_{2,1}$, \ldots, 
$T_{2,p_2-d+1}$, \ldots, $T_{v,1}$, \ldots, $T_{v,p_v-d+1}$ is a shelling of 
${\cal T}(P^{d,k,n})$.
\end{proof}

\begin{corollary}
Let $n\ge k\ge d=2m+1\ge 5$.
Let $\cup[G_j,F_j]$ be the partition of the face lattice of $P^{d,k,n}$ from
the colex shelling, and let $\cup [U_{j,\ell},T_{j,\ell}]$ be the partition
of the face lattice of ${\cal T}(P^{d,k,n})$ from the shelling of 
Theorem~\ref{triang-shell}. Then
\begin{enumerate}
\item $h(P^{d,k,n},x)\ge h'(P^{d,k,n},x)$.
\item The contribution to $h_i(P^{d,k,n})-h'_i(P^{d,k,n})$ from the interval
      $[G_j,F_j]$ is
   $$a_{j,i}=|\{\ell:\mbox{$|U_{j,\ell}|=i$, $1\le \ell\le p_\ell-d$}\}|\ge 0.$$
\end{enumerate}

\end{corollary}

\begin{proof}
The $h$-vector of ${\cal T} ={\cal T}(P^{d,k,n})$ counts the sets $U_{j,\ell}$ 
of each size.
Among these are all the sets $G_j$ counted by the $h'$-vector of $P^{d,k,n}$.
Thus 
\begin{eqnarray*}
h_i(P^{d,k,n})&=&h_i({\cal T})=|\{(j,\ell):|U_{j,\ell}|=i\}|\\
&\ge & |\{(j,\ell):\mbox{$|U_{j,\ell}|=i$ and $\ell=p_j-d+1$}\}|=
h'_i(P^{d,k,n}).
\end{eqnarray*}

Recall that we write ${\cal G}_j$ for the set of faces of $F_j$ not in 
$\cup_{i<j} F_i$; here ${\cal G}_j$ is the set of faces in $[G_j,F_j]$.
Write also ${\cal TG}_j $ for the set of faces of $\cal T$ that are contained
in $F_j$ but not in $\cup_{i<j}F_i$.
By \cite[Corollary 7]{Bayer-weakly}, since $\cal T$ is a shallow triangulation 
of $\partial P^{d,k,n}$, $g(G,x)=\sum (x-1)^{d-1-\dim\sigma}$, where
the sum is over all faces $\sigma$ of $\cal T$ that are contained in $G$ but
not in any proper subface of $G$.
Thus
\begin{eqnarray*}
h({\cal G}_j,x)&=&\sum_{G\in [G_j,F_j]} g(G,x)(x-1)^{d-1-\dim G}\\
&=& \sum_{\sigma\in{\cal TG}_j} (x-1)^{d-1-\dim\sigma}
=\sum_{\ell=1}^{p_\ell-d+1} x^{d-|U_{j,\ell}|}
\end{eqnarray*}
Since $h'({\cal G}_j,x)=x^{d-|G_j|}=x^{d-|U_{j,p_j-d+1}|}$, 
$$\sum_i a_{j,i}x^i=h({\cal G}_j,x)-h'({\cal G}_j,x)=\sum_{\ell=1}^{p_\ell-d}
x^{d-|U_{j,\ell}|},$$
or 
$$a_{j,i}=|\{\ell:\mbox{$|U_{j,\ell}|=i$, $1\le \ell\le p_\ell-d$}\}|\ge 0.$$
\end{proof}

\begin{table}[hbt]
\begin{center}
\begin{tabular}{|l|l|r||l|l|r|}
\hline
$(j,\ell)$ & \multicolumn{1}{c|}{$T_{j,\ell}$} & 
\multicolumn{1}{c||}{$U_{j,\ell}$} & $(j,\ell)$ & 
\multicolumn{1}{c|}{$T_{j,\ell}$} & \multicolumn{1}{c|}{$U_{j,\ell}$}\\
\hline
${1,1}$ & \verb[01234    [ & {\tt $\emptyset$} &
${11,1}$ & \verb[ 1234  7 [ & \verb[27[ \\
${2,1}$ & \verb[012 45   [ & \verb[5[ &
${11,2}$ & \verb[  234  78[ & \verb[78[ \\
${3,1}$ & \verb[0 2345   [ & \verb[35[ &
${12,1}$ & \verb[ 12 45 7 [ & \verb[257[\\
${4,1}$ & \verb[0 23 56  [ & \verb[6[ &
${12,2}$ & \verb[  2 45 78[ & \verb[578[\\
${5,1}$ & \verb[0  3456  [ & \verb[46[ & 
${13,1}$ & \verb[0123  6  [ & \verb[126[\\
${6,1}$ & \verb[01 34 6  [ & \verb[16[ &
${13,2}$ & \verb[ 123  67 [ & \verb[267[\\
${6,2}$ & \verb[ 1 34 67 [ & \verb[7[ &
${13,3}$ & \verb[  23  678[ & \verb[678[\\
${7,1}$ & \verb[01  456  [ & \verb[156[ &
${14,1}$ & \verb[   34 678[ & \verb[4678[\\
${7,2}$ & \verb[ 1  4567 [ & \verb[57[ &
${15,1}$ & \verb[012  56  [ & \verb[1256[\\
${8,1}$ & \verb[  2345  8[ & \verb[8[ &
${15,2}$ & \verb[ 12  567 [ & \verb[2567[\\
${9,1}$ & \verb[  23 56 8[ & \verb[68[ &
${15,3}$ & \verb[  2  5678[ & \verb[5678[\\
${10,1}$ & \verb[   3456 8[ & \verb[468[ &
${16,1}$ & \verb[    45678[ & \verb[45678[\\
\hline
\end{tabular}
\end{center}
\caption{Shelling of triangulation of $P^{5,6,8}$}\label{TP568}
\end{table}

{\bf Example.}
Table \ref{TP568} gives the shelling of the triangulation of $P^{5,6,8}$.
(Refer back to Table~\ref{P568} for the shelling of $P^{5,6,8}$ itself.)
Among the rows $(6,1)$, $(7,1)$, $(11,1)$, $(12,1)$, $(13,1)$,
$(13,2)$, $(15,1)$, $(15,2)$ (rows $(j,\ell)$ that are not
the last row for that $j$), count the $U_{j,\ell}$ of
cardinality $i$ to get $h_i(P^{5,6,8})-h'_i(P^{5,6,8})$.
Note that $U_{13,3}=G_{13}$ (from Table~\ref{P568}), and that
$U_{13,2}=U_{13,3}\setminus\{8\}\cup\{2,7\}$.
The ridges in $T_{13,2}$ are 1236, 1237, 1267, 1367, and 2367.
The first ridge, 1236, falls under Case~1 of the (later part of the) proof of 
Theorem~\ref{triang-shell}; it is contained in the previous facet,
$T_{13,1}$.
The next ridge, 1237, falls under Case~3; it is contained in the
ridge 12378 of $P^{5,6,8}$ in $F_{13}=0123678$, and 12378 also 
contains the ridge 2378 in $T_{13,3}$.
The induction assumption says that 2378 is contained in an earlier facet,
in this case $T_{11,2}$, and 12378 is contained in $F_{11}$.  
Finally, the ridge 1237 is contained in the simplex $T_{11,1}$, part of
the triangulation of $F_{11}$.
The last ridge of $T_{13,2}$ not containing 267 is 1367.
It falls under Case~4.
The set 1367 is contained in the ridge 01367 of $P^{5,6,8}$, contained
in $F_{13}$.  
This ridge is also contained in the earlier facet $F_6$.
The ridge 1367 of the triangulation is contained in the simplex
$T_{6,2}$.

\begin{theorem}
Let $n\ge d+k-1$.
For $1\le i\le d-1$, $h_i(P^{d,k,n})-h_i(P^{d,k,n-1})$ is the number of 
facets $T_{j,\ell}$ of ${\cal T}(P^{d,k,n})$ such that $\max F_j=n-1$
and $|U_{j,\ell}|=i$.
For $1\le i\le (d-1)/2$, this is $k-d+i-1 \choose i-1$.
\end{theorem}

\begin{proof}
Refer to Proposition~\ref{recurs-facets} for a 
description of the facets of $P^{d,k,n}$ in terms of those of $P^{d,k,n-1}$.
For $n\ge d+k-1$, for every facet $P^{d,k,n}$ with maximum element $n$, the 
translation $F-1$ is a facet of $P^{d,k,n-1}$.
(For smaller $n$, a facet of $P^{d,k,n}$ may contain both 0 and $n$, in which 
case $\mbox{lsh(F)}$ is a proper subset of $F-1$.)
The same holds for the simplices $T_{j,\ell}$ triangulating these facets,
and for the sets $U_{j,\ell}$.
The facets of $P^{d,k,n}$ with maximum element at most $n-2$ are facets of 
$P^{d,k,n-1}$,
and the same holds for the corresponding $T_{j,\ell}$ and $U_{j,\ell}$.
The contributions to $h(P^{d,k,n})$ from facets ending in any element
but $n-1$ thus total $h(P^{d,k,n-1})$.
So for $1\le i\le d-1$, $h_i(P^{d,k,n})-h_i(P^{d,k,n-1})$ is the number of 
facets $T_{j,\ell}$ of ${\cal T}(P^{d,k,n})$ such that $\max F_j=n-1$
and $|U_{j,\ell}|=i$.

Now consider the set $\cal S$ of facets $T_{j,\ell}$ of ${\cal T}(P^{d,k,n})$ 
with $\max F_j=n-1$.
For each $T\in{\cal S}$, $T$ is a set of $d$ elements occurring 
consecutively in some $F_j$ with maximum element $n-1$.
So $T$ can be written as
\begin{equation}\label{T}
T=[b,n-k-1]\cup[n-k+1,c]\cup Y\cup[e,b+k],
\end{equation}
where
\begin{enumerate}
\item $n-k-d+1\le b\le n-k-1$;
\item $n-k\le c\le b+d-1$ and $c-n+k$ is even
      (here $c=n-k$ means $[n-k+1,c]=\emptyset$);
\item $Y$ is a paired subset of $[c+2,e-1]$;
\item $e = b+k-1$ if $n-k-b$ is odd, and $e=b+k$ if $n-k-b$ is even; and
\item $|T|=d$.
\end{enumerate}
In these terms, the minimum new face $U$ when $T$ is shelled on is
$U=[b+1,n-k-1]\cup E(Y)\cup \{b+k\}$.

We give a bijection between the facets $T$ in $\cal S$ with $|U|=i$
(where $1\le i\le (d-1)/2$) and the $(k-d)$-element subsets of $[1,k-d+i-1]$.
Let $T$ be as in Equation~\ref{T}.
Then $i=|U|=n-k-b+|Y|/2$.
Write $[c+1,e-1]\setminus Y=\{x_1,x_2,\ldots, x_{k-d}\}$, with the $x_\ell$s
increasing.
(This set has $k-d$ elements because $d=(c-b)+|Y|+(b+k-e+1)$, so
$|[c+1,e-1]\setminus Y|=e-c-1-|Y|=k-d$.)
For each $\ell$, let $y(x_\ell)$ be the number of pairs in $Y$
with both elements less than $x_\ell$.
Let $a_1=n-k-b=i-|Y|/2$.
Set $$A(T)=\{a_1+y(x_\ell)+\ell-1: 1\le \ell\le k-d\}.$$
To see that this is a subset of $[1,k-d+i-1]$, note that the elements of
$A(T)$ form an increasing sequence with minimum element $a_1$ and 
maximum element
$a_1+y(x_{k-d})+(k-d-1) \le a_1+|Y|/2+(k-d-1)=k-d+i-1$.

For the inverse of this map, write a $(k-d)$-element subset of $[1,k-d+i-1]$
as $A =\{a_1,a_2,\ldots, a_{k-d}\}$, with the $a_\ell$s increasing.
Then $1\le a_1\le i$.
Let $$x_1=n-k+d-2i+a_1-\chi(\mbox{$a_1$ odd}).$$
Set 
\begin{eqnarray*}
T(A)&=&[n-k-a_1,n-k-1]\cup[n-k+1,x_1-1]\\
& & {}\cup Y \cup [n-a_1-\chi(\mbox{$a_1$ odd}),n-a_1],
\end{eqnarray*}
where
$$Y= [x_1,n-a_1-1-\chi(\mbox{$a_1$ odd})]\setminus\{x_1+2(a_\ell-a_1)-(\ell-1):
1\le \ell\le k-d\}.$$
We check that this gives a set of the required form.

\noindent(1) Since $1\le a_1\le i\le d-1$,
$n-k-d+1\le n-k-a_1\le n-k-1$.

\noindent(2) 
$x_1-1-n+k=d-2i-1+(a_1-\chi(\mbox{$a_1$ odd}))$, which is nonnegative and even;
$x_1-1=(n-k-a_1+d-1)-(2i-2a_1+\chi(\mbox{$a_1$ odd}))\le n-k-a_1+d-1$.

\noindent(3) 
$Y$ is clearly a subset of $[x_1+1,n-a_1-\chi(\mbox{$a_1$ odd})-1]$.
To see that $Y$ is paired, note that the difference between two consecutive
elements in the removed set is 
$(x_1+2(a_{\ell+1}-a_1)-\ell)-(x_1 +2(a_{\ell}-a_1)-(\ell-1))=2(a_{\ell+1}
-a_\ell)-1$.

\noindent(4) This condition holds by definition.

\noindent(5) Since $A$ is a subset of $[1,k-d+i-1]$,
\begin{eqnarray*}
\lefteqn{x_1+2(a_{k-d}-a_1)-(k-d-1)}\\
&\le & x_1+2(k-d+i-1)-2a_1-(k-d-1)\\
&=&x_1+k-d+2i-2a_1-1=n-a_1-\chi(\mbox{$a_1$ odd})-1.
\end{eqnarray*}
So 
$$\{x_1+2(a_\ell-a_1)-(\ell-1):1\le \ell\le k-d\}\subseteq
[x_1,n-a_1-1-\chi(\mbox{$a_1$ odd})],$$ and
$$|Y|=(n-a_1-\chi(\mbox{$a_1$ odd})-x_1)-(k-d)=2i-2a_1.$$
So $|T(A)|=x_1-(n-k-a_1)+|Y|+ \chi(\mbox{$a_1$ odd})=d$.

Also, in this case $U=[n-k-a_1+1,n-k-1]\cup E(Y)\cup\{n-a_1\}$,
so $|U|=i$.

It is straightforward to check that these maps are inverses.
The main point is that, if $a_\ell=a_1+y(x_\ell)+\ell-1$, then
\begin{eqnarray*}
x_1+2(a_\ell-a_1)-(\ell-1)
&=&x_1+2(y(x_\ell)+\ell-1)-(\ell-1)\\
&=&x_1+2y(x_\ell)+\ell-1=x_\ell.
\end{eqnarray*}
\end{proof}

{\bf Example.}
Consider the ordinary polytope $P^{7,9,15}$.
There are six facets with maximum vertex 14; they are (with sets $G_j$
underlined)
$\{4,5,7,8,9,10,13,\underline{14}\}$,
$\{4,5,7,8,10,\underline{11},13,\underline{14}\}$,
$\{4,5,8,\underline{9},10,\underline{11},13,\underline{14}\}$,
$\{2,3,4,5,7,8,11,\underline{12},13,\underline{14}\}$,
$\{2,3,4,5,8,\underline{9},11,\underline{12},13,\underline{14}\}$,
and
$\{0,1,2,3,4,5,9,\underline{10},11,\underline{12},13,\underline{14}\}$.
Among the 6-simplices occurring in the triangulation of these facets,
six have $|U_{j,\ell}|=3$.
Table~\ref{bijection} gives the bijection from this set of simplices
to the 2-element subsets of $[1,4]$.

\begin{table}[hbt]
\begin{center}
\begin{tabular}{|r|c|c|c|r|c|c|c|c|}
\hline
\multicolumn{1}{|c|}{$T_{j,\ell}$} & $b$ & $c$ & $e$ & \multicolumn{1}{c|}{$Y$} 
& $a_1$ & $x_1,x_2$ & $y(x_i)$ & $A(T_{j,\ell})$\\
\hline
$4,\underline{5},7,8,10,\underline{11},\underline{13}$ 
& 4 & 8 & 13 & $10,11$ & 2 & $9,12$ & $0,1$ & $\{2,4\}$\\
$5,8,\underline{9},10,\underline{11},13,\underline{14}$ 
& 5 & 6 & 13 & $8,9,10,11$ & 1 & $7,12$ & $0,2$ & $\{1,4\}$\\
$3,\underline{4},\underline{5},7,8,11,\underline{12}$
& 3 & 8 & 11 & $\emptyset$ & 3 & $9,10$ & $0,0$ & $\{3,4\}$\\
$4,\underline{5},7,8,11,\underline{12},\underline{13}$
& 4 & 8 & 13 & $11,12$ & 2 & $9,10$ & $0,0$ & $\{2,3\}$\\
$5,8,\underline{9},11,\underline{12},13,\underline{14}$
& 5 & 6 & 13 & $8,9,11,12$ & 1 & $7,10$ & $0,1$ & $\{1,3\}$\\
$5,9,\underline{10},11,\underline{12},13,\underline{14}$
& 5 & 6 & 13 & $9,10,11,12$ & 1 & $7,8$ & $0,0$ & $\{1,2\}$\\
\hline
\end{tabular}
\end{center}
\caption{Bijection with 2-element subsets of $\{1,2,3,4\}$}\label{bijection}
\end{table}

Again, the results of this section hold for even-dimensional multiplexes
as well.

\section{Afterword}
The story of the combinatorics of simplicial polytopes is a beautiful one.
There one finds an intricate interplay among the face lattice of the 
polytope, shellings, the Stanley-Reisner ring and the toric variety,
tied together with the $h$-vector.  
(See, for example, \cite[Lecture 8]{Ziegler}.)
The cyclic polytopes play a special role, serving as the extreme 
examples, and providing the environment in which to build representative
polytopes for each $h$-vector (the Billera-Lee construction \cite{billera-lee}).
In the general case of arbitrary convex polytopes, the various puzzle
pieces have not interlocked as well.
In this paper we made progress on putting the puzzle together for the
special class of ordinary polytopes.
Since the ordinary polytopes generalize the cyclic polytopes, 
a natural next step would be to mimic the Billera-Lee construction, or
Kalai's extension of it \cite{kalai-many}, on
the ordinary polytopes, as a way of generating multiplicial flag vectors.
It would also be interesting to see if there is a ring associated with 
these polytopes, particularly one having a quotient with Hilbert
function equal to the $h'$-polynomial.
Another open problem is to determine the best even-dimensional
analogues of the ordinary polytopes.
They may come from taking vertex figures of odd-dimensional ordinary 
polytopes, or from generalizing Dinh's combinatorial description of
the facets of ordinary polytopes.
Looking beyond ordinary and multiplicial polytopes, we should ask 
what other classes of polytopes have shellings with special
properties that relate to the $h$-vector?

\section*{Acknowledgments}
My thanks go to the folks at University of Washington, the Discrete and 
Computational Geometry program at MSRI and the Diskrete Geometrie group at
TU-Berlin,
who listened to me when it was all speculation.
Particular thanks go to Carl Lee for helpful discussions.


\begin{thebibliography}{10}
\bibitem{Barthel}
Gottfried Barthel, Jean-Paul Brasselet, Karl-Heinz Fieseler, and Ludger Kaup,
\newblock Combinatorial intersection cohomology for fans,
\newblock {\em Tohoku Math. J. (2)}, 54(1):1--41, 2002.

\bibitem{Bayer-weakly} 
Margaret~M. Bayer,
\newblock Equidecomposable and weakly neighborly polytopes,
\newblock {\em Israel J. Math.}, {\bf 81} (1993), 301--320. 

\bibitem{bayer-ordflags}
Margaret~M. Bayer,
\newblock Flag vectors of multiplicial polytopes,
\newblock Preprint 2003-016, Mathematical Sciences Research Institute,
  Berkeley, 2003.

\bibitem{Bayer-Bru-Ste}
Margaret M. Bayer, Aaron M. Bruening, and Joshua D. Stewart,
\newblock A combinatorial study of multiplexes and ordinary polytopes,
\newblock {\em Discrete Comput.\ Geom.}, {\bf 27} (2002), no.~1, 49--63

\bibitem{billera-lee}
Louis~J. Billera and Carl~W. Lee,
\newblock A proof of the sufficiency of {M}c{M}ullen's conditions for
$f$-vectors of simplicial polytopes,
\newblock {\em J.\ Combin.\ Theory Ser.\ A} {\bf 31} (1981), 237--255.

\bibitem{Bisz-mult} 
T.~Bisztriczky,
\newblock On a class of generalized simplices,
\newblock {\em Mathematika}, {\bf 43} (1996), 274--285.

\bibitem{Bisz}
T.~Bisztriczky,
\newblock Ordinary $(2m+1)$-polytopes,
\newblock {\em Israel J. Math.}, 102 (1997), 101--123.

\bibitem{braden}
Tom Braden,
\newblock $g$- and $h$-polynomials of non-rational polytopes---recent progress,
\newblock Abstract at Meeting on Topological and Geometric Combinatorics,
\newblock Mathematisches Forschungsinstitut Oberwolfach,
\newblock April 6--12, 2003.

\bibitem{Dinh}
Thi~Ngoc Dinh,
\newblock {\em Ordinary Polytopes},
\newblock PhD thesis, The University of Calgary, 1999.

\bibitem{kalai-many}
Gil Kalai,
\newblock Many triangulated spheres,
\newblock {\em Discrete Comput.\ Geom.}, {\bf 3} (1988), 1--14.

\bibitem{Karu}
Kalle Karu,
\newblock Hard Lefschetz Theorem for nonrational polytopes,
\newblock arXiv:math.AG/0112087, 2001.

\bibitem{McMullen}
P.~McMullen and G. C. Shephard,
\newblock {\em Convex polytopes and the upper bound conjecture}, 
\newblock Cambridge Univ. Press, London, 1971.

\bibitem{sta80}
Richard~P.\ Stanley,
\newblock The number of faces of simplicial convex polytopes,
\newblock {\em Adv.\ Math.} {\bf 35} (1980), 236--238.

\bibitem{sta87}
Richard~P. Stanley,
\newblock {G}eneralized ${H}$-vectors, intersection cohomology of toric
  varieties, and related results,
\newblock In {\em {C}ommutative algebra and combinatorics ({K}yoto, 1985)},
  volume~11 of {\em Adv. Stud. Pure Math.}, pages 187--213, Amsterdam-New York,
  1987, North-Holland.

\bibitem{sta92}
Richard~P. Stanley,
\newblock Subdivisions and local $h$-vectors,
\newblock {\em J. Amer. Math. Soc.}, {\bf 5} (1992), 805--851.

\bibitem{Ziegler}
G\"{u}nter M. Ziegler,
\newblock {\em Lectures on Polytopes},
\newblock Springer-Verlag, 1995.

\end{thebibliography}
\end{document}